\input ./style/arxiv-general.cfg
\documentclass[aos,MSNbibl,dvips]{arximspdf}
\makeatletter
   \@ifpackageloaded{graphicx}{}{\usepackage{graphicx}}
\makeatother

\doi{10.1214/15-AOS1324}
\volume{43}
\issue{4}
\pubyear{2015}
\firstpage{1774}
\lastpage{1800}
\docsubty{FLA}

\makeatletter
\def\cal{\mathcal}

\newcommand{\eqref}[1]{(\ref{#1})}

\newtheorem{theorem}{Theorem}[section]

\newtheorem{lemma}[theorem]{Lemma}
\newtheorem{corollary}[theorem]{Corollary}
\newproclaim{remark}{Remark}[section]
\newproclaim{example}[theorem]{Example}

\newcommand{\E}{{\mathbb E}}
\newcommand{\R}{{\mathbb R}}
\renewcommand{\P}{{\mathbb P}}
\newcommand{\C}{{\mathcal{C}}}
\newcommand{\nt}{{\mathbb N}}
\newcommand{\mix}{{\mathfrak{R}}}
\newcommand{\nbr}{{\mathfrak{N}}}
\newcommand{\M}{{\mathcal{M}}}
\newcommand{\K}{{\mathcal{K}}}
\newcommand{\Ps}{{\mathcal{P}}}
\newcommand{\F}{{\cal F}}
\newcommand{\lin}{{\mathcal{L}}}
\newcommand{\ham}{{\Upsilon}}

\def\argmin{\mathop{\operatorname{argmin}}}
\def\qt#1{\mbox{#1}}
\makeatother

\begin{document}
\begin{frontmatter}

\title{On risk bounds in isotonic and other shape restricted
regression problems}
\runtitle{On risk bounds in shape restricted regression}

\begin{aug}
\author[A]{\fnms{Sabyasachi}~\snm{Chatterjee}\ead[label=e1]{sabyasachi@galton.uchicago.edu}},
\author[B]{\fnms{Adityanand}~\snm{Guntuboyina}\corref{}\thanksref{T2}\ead[label=e2]{aditya@stat.berkeley.edu}}\\
\and
\author[C]{\fnms{Bodhisattva}~\snm{Sen}\thanksref{T3}\ead[label=e3]{bodhi@stat.columbia.edu}}
\runauthor{S. Chatterjee, A. Guntuboyina and B. Sen}

\thankstext{T2}{Supported by NSF Grant DMS-13-09356.}
\thankstext{T3}{Supported by NSF Grants DMS-11-50435 and AST-1107373.}
\affiliation{University of Chicago, University of California,
Berkeley and Columbia University}

\address[A]{S. Chatterjee\\
Department of Statistics\\
University of Chicago\\
5734 S.~University Avenue \\
Chicago, Illinois 60637\\
USA\\
\printead{e1}}

\address[B]{A. Guntuboyina\\
Department of Statistics\\
University of California, Berkeley\\
423 Evans Hall\\
Berkeley, California 94720\\
USA\\
\printead{e2}\\
\phantom{E-mail:\ }}

\address[C]{B. Sen\\
Department of Statistics\\
Columbia University\\
1255 Amsterdam Avenue\\
New York, New York 10027\\
USA\\
\printead{e3}}
\end{aug}

%
\received{\smonth{5} \syear{2014}}
%
\revised{\smonth{2} \syear{2015}}

%
\begin{abstract}
We consider the problem of estimating an unknown $\theta\in\R^n$
from noisy observations under the constraint that $\theta$ belongs to
certain convex polyhedral cones in $\R^n$. Under this setting, we
prove bounds for the risk of the least squares estimator (LSE). The
obtained risk bound behaves differently depending on the true sequence
$\theta$ which highlights the adaptive behavior of~$\theta$. As
special cases of our general result, we derive risk bounds for the LSE
in univariate isotonic and convex regression. We study the
risk bound in isotonic regression in greater detail: we show that
the isotonic LSE converges at a whole range of rates from $\log n/n$
(when $\theta$ is constant) to $n^{-2/3}$ (when $\theta$ is
\textit{uniformly increasing} in a certain sense). We argue that the
bound presents a benchmark for the risk of any estimator in isotonic
regression by proving nonasymptotic local minimax lower bounds. We
prove an analogue of our bound for model misspecification where the
true $\theta$ is not necessarily nondecreasing.
\end{abstract}

%
\end{frontmatter}

\section{Introduction}\label{sec1}
Shape constrained regression involves estimating a vector $\theta=
(\theta_1, \ldots, \theta_n) \in\R^n$ from observations
%
\begin{equation}
\label{eq:Mdl2} Y_i = \theta_i + \varepsilon_i\qquad
\qt{for $i = 1, \ldots, n$},
\end{equation}
where $\theta$ lies in a known convex polyhedral cone
$\K\subseteq\R^n$ and $\varepsilon_1, \ldots, \varepsilon_n$ are
i.i.d. mean zero errors with finite variance. It may be recalled that convex
polyhedral cones are sets of the form
%
\begin{equation}
\label{bla} \K:= \bigl\{\theta\in\R^n\dvtx A \theta\geq0 \bigr\},
\end{equation}
where $A$ is a matrix of order $m \times n$ and $\alpha= (\alpha_1,
\ldots, \alpha_m) \geq0$ means that $\alpha_i \geq0$ for each
$i$. Basic properties of convex polyhedral
cones can be found, for example, in \cite{Schrijver}, Chapters 7 and 8.

In this paper, we focus on such problems when the cone $\K$ is of the
special form
%
\begin{equation}
\label{gv} \K^n_{r, s} := \Biggl\{\theta\in
\R^n \dvtx \sum_{j = -r}^s
w_j \theta_{t + j} \geq0 \mbox{ for all } 1 + r \leq t \leq n
- s \Biggr\},
\end{equation}
for some known integers $r \geq0$ and $s \geq1$ and nonnegative
weights $w_{j}, -r \leq j \leq s$. Here the integers $r$ and $s$ and
the weights $w_j, -r \leq j \leq s$ do not depend on $n$. Note that
when $n < 1 + r + s$, the condition in the definition of $\K_{r,
s}^n$ is vacuous so that $\K_{r, s}^n = \R^n$. The dependence of the
cone on the weights $\{w_j\}$ is suppressed in the notation $\K_{r,
s}^n$.

The following shape constrained regression problems are special
instances of our general setup:
\begin{longlist}[(1)]
\item[(1)] When $r = 0, s = 1, w_0 = -1$ and $w_1 = 1$, the cone
in \eqref{gv} consists of all nondecreasing
sequences
\[
\M:= \bigl\{\theta\in\R^n \dvtx \theta_1 \leq\cdots\leq
\theta_n\bigr\}.
\]
Estimation problem \eqref{eq:Mdl2} then becomes the well-known
isotonic regression problem.

\item[(2)] When $r = 1, s = 1$, $w_{-1} = w_1 = 1$ and $w_0 = -2$, the cone
in \eqref{gv} consists of all convex sequences $\C:= \{\theta\in \R
^n\dvtx 2 \theta_i \leq\theta_{i-1} + \theta_{i+1}, i = 2, \ldots, n-1\}
$. Then \eqref{eq:Mdl2} reduces to the usual convex regression problem
with equally
spaced design points.
\item[(3)] $k$-monotone regression corresponds to $\K:= \{\theta\in
\R^n\dvtx \nabla^k \theta\geq0\}$ where $\nabla\dvtx \R^n
\rightarrow\R^n$ is given by $\nabla(\theta) := (\theta_2 -
\theta_1, \theta_3 - \theta_2, \ldots, \theta_n - \theta_{n-1},
0)$, and $\nabla^k$ represents the $k$-times composition of $\nabla$.
This is also a special case of \eqref{gv}.
\end{longlist}

Our object of interest in this paper is the least squares estimator
(LSE) for $\theta$ under the constraint $\theta\in\K$. It is given by
$\hat{\theta}(Y; \K)$ where $Y = (Y_1, \ldots, Y_n)$ is the
observation vector and
%
\begin{equation}
\label{eq:LSE} \hat{\theta}(y; \K) := \argmin_{\theta\in\K} \|\theta- y
\|^2 \qquad\qt{for $y \in\R^n$},
\end{equation}
where $\|\cdot\|$ denotes the usual Euclidean norm in $\R^n$. A
natural measure of how well $\hat{\theta}(Y; \K)$ estimates $\theta$
is $\ell^2(\theta, \hat{\theta}(Y; \K))$ where
%
\begin{equation}
\label{eq:loss} \ell^2(\alpha, \beta) := \frac{1}{n} \| \alpha-
\beta\|^2 = \frac{1}{n} \sum_{i=1}^n
(\alpha_i - \beta_i)^2
\end{equation}
with $\alpha= (\alpha_1, \ldots, \alpha_n)$ and $\beta= (\beta_1,
\ldots, \beta_n)$. As $\ell^2(\theta, \hat{\theta}(Y; \K))$ is random
we study its expectation $\E_{\theta} \ell^2(\theta, \hat{\theta
}(Y;\K))$ which is referred to as the risk of the estimator $\hat
{\theta}(Y; \K)$.

This paper has two aims: (1) For every cone $\K^n_{r, s}$, we prove
upper bounds for the risk $\E_{\theta} \ell^2(\theta, \hat{\theta
}(Y; \K^n_{r,s}))$ as $\theta$ varies in $\K_{r, s}^n$; (2) we
isolate the risk bound for the special case of isotonic regression
(when $\K_{r,s}^n = \M$) and study its properties in more detail.

\subsection{Upper bounds on the risk of \texorpdfstring{$\hat{\theta}(Y; \K^n_{r, s})$}
{hat{theta}(Ymathcal{K}n{r,s})}}
The first part of the paper will be about bounds for the risk
$\E_{\theta} \ell^2(\theta, \hat{\theta}(Y; \K_{r, s}^n))$. Our
bounds will involve the \textit{statistical dimension} of the cone $\K_{r,
s}^n$. For a cone $\K$, defined as in \eqref{bla}, its statistical
dimension is given by
%
\begin{equation}
\label{sd1} \delta(\K) := \E D(Z; \K) \qquad\qt{where }   D(y; \K) := \sum
_{i=1}^n \frac{\partial}{\partial y_i} \hat{
\theta}_i(y; \K)
\end{equation}
and $Z = (Z_1, \ldots, Z_n)$ is a vector whose components are
independent standard normal random variables. Note that the quantity
$D(y; \K)$ is well defined because $\hat{\theta}(y; \K)$ is a
1-Lipschitz function of $y$; see \cite{MW00}. It was argued in
\cite{MW00} that
$D(Y;\K)$ provides a measure of the effective dimension of the
model. For example, if $\K$ is a linear space of dimension $d$, then
$\hat\theta(y;\K) = Q Y$, where $Q$ is the projection matrix onto
$\K$, and $D(y;\K) = \operatorname{trace}(Q) = d$ for all $y$. It was also
shown in \cite{MW00} that $D(Y;\K)$ is the number of distinct values
among $\hat\theta_1, \ldots, \hat\theta_n$ for isotonic
regression. The term statistical dimension for $\delta(\K)$ was first
used in \cite{LivEdge}; however, the definition of $\delta(\K)$
in \cite{LivEdge} is different from \eqref{sd1}. For connections
between the two definitions and more discussion on the notion of
statistical dimension, see Section~\ref{genres}.

We are now ready to describe our main result which bounds $\E_{\theta
} \ell^2(\theta, \hat{\theta}(Y;\break  \K_{r, s}^n))$ for $\theta\in\K
_{r, s}^n$. For each $\theta\in\K_{r, s}^n$, let $k(\theta)$ denote the
number of inequalities among $\sum_{j=-r}^s w_j \theta_{t + j} \geq0$
for $1+r \leq t \leq n-s$ that are strict. In Theorem~\ref{zen}, we
prove that for every $\theta\in\K_{r, s}^n$,
%
\begin{equation}
\label{zen.eq} \E_{\theta} \ell^2\bigl(\theta, \hat{\theta}
\bigl(Y; \K^n_{r, s}\bigr)\bigr) \leq6 \inf
_{\alpha\in\K_{r, s}^n} \biggl(\ell^2(\theta, \alpha) +
\frac{\sigma^2 (1 + k(\alpha))}{n} \delta\bigl(\K_{r, s}^n\bigr) \biggr)
\end{equation}
under the assumption that $\varepsilon_1, \ldots, \varepsilon_n$ are
independent normally distributed random variables with mean zero and
variance $\sigma^2$. This bound behaves differently depending on the
form of the true sequence $\theta$ and thus describes the adaptive
behavior of the LSE; for more details on the inequality, see
Section~\ref{genres}. The proof of Theorem~\ref{zen} uses the
characterization properties of the projection operator on a closed
convex cone. We prove a series of auxiliary results leading to the
proof of Theorem~\ref{zen}; these results hold for any polyhedral cone
$\K$ (not necessarily of the form $\K^{n}_{r,s}$) and are of
independent interest.

\subsection{On risk bounds in isotonic regression} The second part of
the paper is exclusively on isotonic
regression. Even in this special case, inequality \eqref{zen.eq}
appears to be new. We provide a reformulation of \eqref{zen.eq} that
bounds the risk of $\hat{\theta}(Y; \M)$ using the variation of
$\theta$ across subsets of $\{1, \ldots, n\}$. This results in an
inequality that is more interpretable and makes comparison with
previous inequalities in isotonic regression more transparent.

To state this bound, we need some notation. Specializing the notation
$k(\theta)$ for $\theta\in\K^n_{r, s}$ to the cone $\M$, we get
$k(\theta)$ equals the number of inequalities $\theta_i \leq
\theta_{i+1}$ for $i = 1, \ldots, n-1$ that are strict. By an
abuse of notation, we extend this notation to \textit{interval
partitions} of $n$. An interval partition $\pi$ of $n$ is a finite
sequence of positive integers that sum to $n$. In combinatorics this
is called a composition of $n$. Let the set of all interval
partitions $\pi$ of $n$ be denoted by $\Pi$. Formally, $\Pi$ can be
written as
\[
\Pi:= \Biggl\{(n_1,n_2,\ldots,n_{k+1})\dvtx k
\geq0, n_i \in\nt\mbox{ and } \sum_{i=1}^{k+1}
n_i = n \Biggr\}.
\]
For each $\pi= (n_1, \ldots, n_{k+1})$ $\in\Pi$, let $k(\pi) := k$.

For every $\theta= (\theta_1, \ldots, \theta_n) \in\M$, there exist
integers $k \geq0$ and $n_1, \ldots, n_{k+1} \geq1$ with $n_1 +
\cdots
+ n_{k+1} = n$ such that $\theta$ is constant on each set $\{j\dvtx s_{i-1} + 1 \leq j \leq s_i \}$ for $i = 1, \ldots, k+1$, where $s_0 :=
0$ and $s_i = n_1 + \cdots+ n_i$. We refer to this interval partition
$\pi_{\theta} := (n_1, \ldots, n_{k+1})$ as the interval partition
\textit{generated} by $\theta$. Note that $k(\pi_{\theta})$ precisely
equals $k(\theta)$, the number of inequalities $\theta_i \leq
\theta_{i+1}$, for $i = 1, \ldots, n-1$, that are strict.

For every $\theta\in\M$ and $\pi:= (n_1, \ldots, n_{k+1}) \in\Pi
$, we
define
\[
V_{\pi}(\theta) = \max_{1 \leq i \leq k+1} ( \theta_{s_i}
- \theta_{s_{i-1} + 1} ), %
\]
where $s_0 := 0$ and
$s_i = n_1 + \cdots+ n_i$ for $1 \leq i \leq k+1$. $V_{\pi}(\theta)$ can
be treated as measure of variation of $\theta$ with respect to the
partition $\pi$. An important property is that
$V_{\pi_{\theta}}(\theta) = 0$ for every $\theta\in\M$. For the
trivial partition $\pi= (n)$, it is easy to see that $k(\pi) = 0$ and
$V_{\pi}(\theta) = V(\theta) = \theta_n - \theta_1$.

With this notation, our main result for isotonic regression states
that
%
\begin{equation}
\label{pyaa} \E_{\theta} \ell^2 \bigl( \theta, \hat{\theta}(Y;
\M) \bigr) \leq R(n; \theta),
\end{equation}
where
\[
R(n; \theta) = 4 \inf_{\pi\in\Pi} \biggl( V^2_{\pi}(
\theta) + \frac{4 \sigma^2 (1 + k(\pi))}{n} \log\frac{en}{1 + k(\pi)} \biggr).
\]
This inequality is very similar to \eqref{zen.eq}; see
Remark~\ref{remp} for the connections. The LSE, $\hat{\theta}(Y; \M
)$, in isotonic regression has the explicit
formula \eqref{eq:Charac}. This formula is commonly known as the
min--max formula; see \cite{RWD88}, Chapter~1. Using this formula, we
prove inequality \eqref{pyaa} in Section~\ref{Isoreg}. We only use
the fact that $\varepsilon_1, \ldots, \varepsilon_n$ are i.i.d. with mean
zero and variance $\sigma^2$ [normality of $\varepsilon_1, \ldots,
\varepsilon_n$ is not needed here unlike inequality \eqref{zen.eq} for
which we require
normality].

Inequality \eqref{pyaa} appears to be new, even
though there is a huge literature on univariate isotonic
regression. To place this inequality in a proper historical context,
we give a brief overview of existing theoretical results on isotonic
regression in Section~\ref{Isoreg}. The strongest previous bound on
$\E_{\theta} \ell^2(\theta, \hat{\theta}(Y; \M))$ is due
to Zhang~\cite{Zhang02}, Theorem~2.2, who showed that
%
\begin{equation}
\label{motw} \E_{\theta} \ell^2\bigl(\theta, \hat{\theta}(Y;
\M)\bigr) \lesssim R_{Z}(n; \theta),
\end{equation}
where
\[
R_Z(n; \theta) := \biggl(\frac{\sigma^2 V(\theta)}{n} \biggr)^{2/3} +
\frac{\sigma^2 \log n}{n}
\]
with
\[
V(\theta) := \theta_n - \theta_1.
\]
Here, by the symbol
$\lesssim$ we mean $\leq$ up to a multiplicative constant. The
quantity $V(\theta)$ is known as the variation of the sequence
$\theta$.

Our inequality \eqref{pyaa} compares favorably with \eqref{motw} in
certain cases. To see this, suppose, for example, that $\theta_j = I\{
j >
n/2\}$ (here $I$ denotes the indicator function) so that $V(\theta)
= 1$. Then $R_Z(n; \theta)$ is essentially $(\sigma^2/n)^{2/3}$ while
$R(n; \theta)$ is much smaller because it is at most $(32 \sigma^2/n)
\log(en/2)$ as can be seen by taking $\pi= \pi_{\theta}$
in the definition of $R(n; \theta)$ [note that $k(\theta) = 1$].

More generally by taking $\pi= \pi_{\theta}$ in the infimum of the
definition of $R(n; \theta)$, we obtain
%
\begin{equation}
\label{amon} \E_{\theta} \ell^2\bigl(\theta, \hat{\theta}(Y;
\M)\bigr) \leq\frac{16(1 +
k(\theta)) \sigma^2}{n} \log\frac{en}{1 + k(\theta)},
\end{equation}
which is a stronger bound than \eqref{motw} when $k(\theta)$ is
small. The reader may observe that $k(\theta)$ is small precisely when
the differences $\theta_i - \theta_{i-1}$ are
sparse. Inequality \eqref{pyaa} can be stronger than \eqref{motw} even
in situations when $k(\theta)$ is not small; see Remark~\ref{ach} for
an example.

We study properties of $R(n; \theta)$ in Section~\ref{hd}. In
Theorem~\ref{sabya2}, we show that $R(n; \theta)$ is bounded from
above by a multiple of $R_Z(n; \theta)$ that is at most logarithmic in
$n$. This implies that our inequality \eqref{pyaa} is always only
slightly worse off than \eqref{motw} while being much better in the
case of certain sequences $\theta$. We also show in Section~\ref{hd}
that the risk bound $R(n; \theta)$ behaves differently, depending on
the form of the true sequence $\theta$. This means that
bound \eqref{pyaa} demonstrates adaptive behavior of the LSE. One gets
a whole range of rates from $(\log n)/n$ (when $\theta$ is constant)
to $n^{-2/3}$ up to logarithmic factors in the worst case [this worst
case rate corresponds to the situation where $\min_i (\theta_i -
\theta_{i-1}) \gtrsim1/n$]. Bound \eqref{pyaa} therefore presents
a bridge between the two terms in the formula for $R_Z(n; \theta)$.


In addition to being an upper bound for the risk of the LSE, we
believe that the quantity $R(n; \theta)$ also acts as a benchmark for
the risk of any estimator in isotonic regression. By this, we mean
that, in a certain sense, no estimator can have risk that is
significantly better than $R(n; \theta)$. We substantiate this claim
in Section~\ref{lmnop} by proving lower bounds for the \textit{local
minimax risk} near the ``true'' $\theta$. For $\theta\in\M$, the
quantity
\[
\mix_n(\theta) := \inf_{\hat{t}} \sup
_{t \in\nbr(\theta)} \E_t \ell^2(t,\hat{t})
\]
with
\[
\nbr(\theta) := \bigl\{t \in\M\dvtx \ell_{\infty}^2 (t,\theta)
\lesssim R(n; \theta) \bigr\}
\]
will be called the local minimax risk at $\theta$; see
Section~\ref{lmnop} for the rigorous definition of the neighborhood
$\nbr(\theta)$ where the multiplicative constants hidden by the
$\lesssim$ sign are explicitly given. In the above display
$\ell_\infty$ is defined as $\ell_{\infty} (t,\theta) :=
\max_{i} |t_i - \theta_i|$. The infimum here is over all
possible estimators $\hat{t}$. $\mix_n(\theta)$ represents the
smallest possible (supremum) risk under the knowledge that the true
sequence $t$ lies in the neighborhood $\nbr(\theta)$. It provides a
measure of the difficulty of estimation of $\theta$. Note that the
size of the neighborhood $\nbr(\theta)$
changes with $\theta$ (and with~$n$) and also reflects the
difficulty level of the problem.

Under each of the two following setups for $\theta$, and the
assumption of normality of the errors, we show that $\mix_n(\theta)$
is bounded from below by $R(n; \theta)$ up to multiplicative
logarithmic factors of $n$. Specifically:
\begin{longlist}[(1)]
\item[(1)] 
when the increments of $\theta$
(defined as $\theta_i - \theta_{i-1}$, for $i = 2, \ldots, n$) grow
like $1/n$, we prove in Theorem~\ref{oval} that
%
\begin{equation}
\label{inmi1} \mix_n(\theta) \gtrsim \biggl(\frac{\sigma^2 V(\theta)}{n}
\biggr)^{2/3} \gtrsim  \frac{R(n; \theta)}{\log(4n)};
\end{equation}
\item[(2)] 
when $k(\theta) = k$ and the $k$ values of
$\theta$
are sufficiently well-separated, we show in Theorem~\ref{fani} that
%
\begin{equation}
\label{inmi2} \mix_n(\theta) \gtrsim R(n; \theta) \biggl(\log
\frac{en}{k} \biggr)^{-2/3}.
\end{equation}
\end{longlist}
Because $R(n, \theta)$ is an upper bound for the risk of the LSE and also
is a local minimax lower bound in the above sense, our results imply
that the LSE is near-optimal in a local nonasymptotic minimax sense.
Such local minimax bounds are in the spirit of Cator \cite{Cator2011}
and Cai and Low \cite{CaiLowFwork}, who worked with the problems of estimating
monotone and convex functions respectively at a point. The difference
between these works and our own is that we focus on the global estimation
problem. In other words, \cite{Cator2011} and \cite{CaiLowFwork}
prove local minimax bounds for the local (pointwise) estimation
problem while we prove local minimax bounds for the global estimation
problem. 

We also study the performance of the LSE in isotonic regression under
model misspecification when the true sequence $\theta$ is not
necessarily nondecreasing. Here we prove in Theorem~\ref{thm:Misrafi}
that $\E_{\theta} \ell^2(\tilde\theta, \hat{\theta}(Y; \M))
\leq R(n;\tilde \theta)$ where $\tilde\theta$ denotes the
nondecreasing projection of $\theta$; see Section~\ref{Misspec} for
its definition. This should be contrasted with the risk bound of Zhang
\cite{Zhang02} who proved
that $\E_{\theta} \ell^2(\tilde{\theta},\hat{\theta}(Y; \M))
\lesssim
R_Z(n;\tilde \theta)$. As before our risk bound is at most, slightly
worse (by a multiplicative logarithmic factor in $n$) than $R_Z$, but
is much better when $k(\tilde\theta)$ is small. We describe two
situations where $k(\tilde\theta)$ is small: when $\theta$ itself
has few constant blocks [see \eqref{eq:block} and Lemma~\ref{sabya3}]
and when $\theta$ is nonincreasing [in which case $k(\tilde
\theta)=1$; see Lemma~\ref{appr}].

\subsection{Organization of the paper} The paper is organized as
follows: In Section~\ref{genres} we state and prove our main upper
bound for the risk $\E_\theta\ell^2(\theta,\hat\theta(Y; \K
^n_{r, s}))$. In Section~\ref{Isoreg} we give a direct proof of the
risk bound \eqref{pyaa} for isotonic regression without assuming
normality of the errors. We investigate the behavior of $R(n;\theta)$,
the right-hand side of \eqref{pyaa}, for different values of the true sequence
$\theta$ and compare it with $R_Z(n;\theta)$, the right-hand
of \eqref{motw}, in Section~\ref{hd}. Local minimax lower bounds for
isotonic regression are proved in Section~\ref{lmnop}. We study the
performance of the isotonic LSE under model misspecification in
Section~\ref{Misspec}. 
The supplementary material \cite{ChatEtAlSupp14} gives the proofs of
some of the results in the paper.

\section{A general risk bound for the projection on closed convex
polyhedral cones}\label{genres}

The goal of this section is to prove inequality \eqref{zen.eq}. Let us
first review the well-known characterization of the LSE under the
constraint $\theta\in\K$ for an arbitrary convex polyhedral cone
$\K$. This LSE is denoted by $\hat{\theta}(Y; \K)$ and is defined
in \eqref{eq:LSE}. The function $y \mapsto\hat{\theta}(y; \K)$ is
well defined [because for each $y$ and $\K$, the quantity
$\hat{\theta}(y; \K)$ exists uniquely by the Hilbert projection
theorem], nonlinear in $y$ (in general) and can be characterized by
(see, e.g., \cite{Bertsekas03}, Proposition~2.2.1)
%
\begin{equation}\qquad
\label{ctz} \hat{\theta}(y; \K) \in\K ,\qquad \bigl\langle y - \hat{\theta
}(y; \K), \hat{\theta}(y; \K) \bigr\rangle= 0
\quad\mbox{and}\quad \bigl\langle y - \hat{
\theta}(y; \K), \omega\bigr\rangle\leq0
\end{equation}
for all $\omega\in\K$.

Inequality \eqref{zen.eq} involves the notion of statistical
dimension [defined in \eqref{sd1}]. The statistical dimension is an
important summary parameter for cones,
and it has been used in shape-constrained regression \cite{MW00}
and compressed sensing \cite{LivEdge,OmyakHassibi2013}. It is
closely related to the Gaussian width of $\K$, which is an important
quantity in geometric functional analysis (see, e.g., \cite{VershyninGFA},
Chapter~4) and which has also been used to prove
recovery bounds in compressed sensing \cite
{RudelsonVershynin2006,ChandraFOCS,Stojnic,LivEdge,OmyakHassibi2013}.
See \cite{LivEdge}, Section~10.3, for the
precise connection between the statistical dimension and the Gaussian width.

An alternative definition of the
statistical dimension $\delta(\K)$ of an arbitrary convex polyhedral
cone is given by
%
\begin{equation}
\label{sd2} \delta(\K) = \E\bigl\|\hat{\theta}(Z; \K)\bigr \|^2,
\end{equation}
where $Z = (Z_1, \ldots, Z_n)$ is a vector whose components are
independent standard normal random variables.
The equivalence of \eqref{sd1} and \eqref{sd2} was observed
by Meyer and Woodroofe \cite{MW00}, proof of Proposition~2. It is
actually an easy
consequence of Stein's lemma because the second identity
in \eqref{ctz} implies $\E\|\hat{\theta}(Z; \K)\|^2 = \E\langle
Z, \hat{\theta}(Z; \K) \rangle$, and therefore, Stein's lemma on the
right-hand side gives the equivalence of \eqref{sd1} and \eqref{sd2}.

We are now ready to prove our main result, Theorem~\ref{zen}, which
gives inequality~\eqref{zen.eq}. This theorem applies to any cone of
the form \eqref{gv}. For the proof of Theorem~\ref{zen}, we state
certain auxiliary results (Lemmas \ref{lili}, \ref{hol}
and \ref{apps}), whose proofs can be found in the supplementary
material \cite{ChatEtAlSupp14}. These supplementary results hold
for any polyhedral cone \eqref{bla}.

\begin{theorem}\label{zen}
Fix $n \geq1$, $r \geq0$ and $s \geq1$. Consider the problem of
estimating $\theta\in\K_{r, s}^n$ from \eqref{eq:Mdl2} for
independent $N(0, \sigma^2)$ errors $\varepsilon_1, \ldots,
\varepsilon_n$. Then inequality~\eqref{zen.eq} holds for every $\theta
\in\K_{r, s}^n$ with $k(\theta)$ denoting the number of inequalities
among $\sum_{j = -r}^s w_j \theta_{t + j} \geq0$, for $1 + r \leq t
\leq n - s$, that are strict.
\end{theorem}

Before we prove Theorem~\ref{zen} the following remarks are in order.

\begin{remark}[(Stronger version)]
From the proof of Theorem~\ref{zen}, it will be clear that the risk
of the LSE satisfies a stronger inequality
than \eqref{zen.eq}. For $\alpha\in\K_{r, s}^n$ with $k(\alpha) =
k$, let $1+r \leq t_1 < \cdots< t_k \leq n - s$ denote the values of
$t$ for which the inequalities $\sum_{j = -r}^s w_j \alpha_{t + j}
\geq0$ are strict. Let
%
\begin{equation}
\label{tde} \tau(\alpha) := \delta\bigl(\K_{r, s}^{t_1 - 1 + s}\bigr)
+ \delta\bigl(\K_{r,
s}^{t_2 - t_1}\bigr) + \cdots+ \delta\bigl(
\K_{r, s}^{t_k - t_{k-1}}\bigr) + \delta\bigl(\K_{r, s}^{n - t_k - s + 1}
\bigr).
\end{equation}
The proof of Theorem~\ref{zen} will imply that
%
\begin{equation}
\label{zens} \E_{\theta} \ell^2 \bigl(\theta, \hat{\theta}
\bigl(Y; \K^n_{r, s}\bigr) \bigr) \leq6 \inf
_{\alpha\in\K_{r, s}^n} \biggl(\ell^2(\theta, \alpha) +
\frac{\sigma^2}{n} \tau(\alpha) \biggr).
\end{equation}
The observation that $\delta(\K_{r, s}^n)$ is increasing in $n$ (note
that the weights $w_j, -r \leq j \leq s$, do not depend on
$n$) implies that $\tau(\alpha) \leq(1 + k(\alpha)) \delta(\K_{r,
s}^n)$ for all $\alpha\in\K_{r, s}^n$, and hence
inequality \eqref{zens} is stronger than \eqref{zen.eq}.
\end{remark}

\begin{remark}[(Connection to the facial structure of $\K_{r, s}^n$)]
 Every convex polyhedral cone \eqref{bla} has a well-defined facial
structure. Indeed, a standard result (see, e.g.,  \cite{Schrijver},
Section~8.3) states that a subset $F$ of a convex
polyhedral cone $\K$, as defined in \eqref{bla}, is a face if and
only if $F$ is nonempty and $F = \{\theta\in\K\dvtx \tilde{A} \theta=
0\}$ for some $\tilde{m} \times n$ matrix $\tilde{A}$ whose rows are
a subset of the rows of $A$. The dimension of $F$ equals $n - \rho
(\tilde{A})$ where $\rho(\tilde{A})$ denotes the rank of $\tilde
{A}$. It is then clear that if $\theta\in\K_{r, s}^n$ is in a
low-dimensional face of $\K_{r, s}^n$, then $k(\theta)$ must be
small. Now if
$\delta(\K_{r, s}^n)$ is at most logarithmic in $n$ (which is indeed
the case for the case of isotonic and convex regression; see
Examples \ref{ir} and \ref{conre}), then bound \eqref{zen.eq} implies
that the risk of the LSE is bounded from above by the parametric rate
$\sigma^2/n$ (up to multiplicative logarithmic factors in $n$) provided
$\theta$ is in a low-dimensional face of $\K_{r, s}^n$. Therefore,
the LSE
automatically adapts to vectors in low-dimensional faces of $\K^n_{r,
s}$. For general $\theta$, the risk is bounded from above by a
combination of how close $\theta$ is to a $k$-dimensional face of
$\K^n_{r, s}$ and $\sigma^2 \delta(\K_{r, s}^n)(1 + k)/n$ as $k$
varies.
\end{remark}

\begin{example}[(Isotonic regression)]\label{ir}
Isotonic regression corresponds to $r = 0, s = 1, w_0 = -1$ and $w_1
= 1$ so that $\K^n_{r,s}$ becomes $\M$. It turns out that the
statistical dimension of this cone satisfies
%
\begin{equation}
\label{frg} \delta(\M) = 1 + \frac{1}{2} + \cdots+ \frac{1}{n}\qquad
\qt{for every $n \geq1$,}
\end{equation}
which immediately implies that $\delta(\M) \leq\log
(en)$. This can be proved using symmetry arguments formalized in the
theory of finite reflection groups; see \cite{LivEdge}, Appendix
C.4, where the proof of \eqref{frg} is
sketched.

Now let $\alpha= (\alpha_1, \ldots, \alpha_n) \in\M$ with
$k(\alpha) = k$. Then there exist integers $n_1, \ldots, n_{k+1} \geq
1$ with $n_1 + \cdots+ n_{k+1} = n$ such that $\alpha$ is constant on
each set $\{j\dvtx s_{i-1} + 1 \leq j \leq s_i \}$, for $i = 1, \ldots,
k+1$, where $s_0 := 0$ and $s_i = n_1 + \cdots+ n_i$. It is easy to
check then that the quantity $\tau(\alpha)$ defined in \eqref{tde} equals
\[
\tau(\alpha) = \delta\bigl(\M^{n_1}\bigr) + \cdots+ \delta\bigl(
\M^{n_{k+1}}\bigr),
\]
where $\M^i :=  \{(\theta_1, \ldots, \theta_i) \in\R^i\dvtx \theta_1 \leq
\cdots\leq\theta_i  \}$ is the monotone cone in $\R^i$.
Inequality \eqref{frg} then gives
\[
\tau(\alpha) \leq\sum_{i=1}^{k+1}
\log(en_i) \leq(k+1) \log \biggl( \frac{en}{k+1} \biggr)
\]
because of the concavity of $x \mapsto\log x$. Inequality \eqref
{zens} therefore gives
%
\begin{equation}\quad
\label{hh} \E_{\theta} \ell^2 \bigl(\theta, \hat{\theta}(Y;
\M) \bigr) \leq6 \inf_{\alpha\in\M} \biggl(\ell^2(\theta,
\alpha) + \frac{\sigma^2 (k(\alpha)+1)}{n} \log \frac{en}{k(\alpha) + 1} \biggr).
\end{equation}
This inequality is closely connected to \eqref{pyaa}, as we describe
in detail in Remark~\ref{remp}. Note that we require normality of
$\varepsilon_1, \ldots, \varepsilon_n$. In Section~\ref{Isoreg}, we prove
an inequality which gives a variant of inequality \eqref{hh} with
different multiplicative constants but without the assumption of normality.

Inequality \eqref{hh} can be restated in the following way. For each $0
\leq k \leq n-1$, let $\Ps_k$ denote the set of all sequences $\alpha
\in\M$ with $k(\alpha) \leq k$. With this notation,
inequality \eqref{hh} can be rewritten as
%
\begin{equation}\quad
\label{vk} \E_{\theta} \ell^2\bigl(\theta, \hat{\theta}(Y;
\M)\bigr) \leq6 \min_{0
\leq k \leq n-1} \biggl[\inf_{\alpha\in\Ps_k}
\ell^2(\theta,\alpha) + \frac
{\sigma^2
(k+1)}{n} \log\frac{en}{k+1}
\biggr].
\end{equation}

Bound \eqref{vk} reflects adaptation of the LSE with respect to
the classes $\Ps_k$. Such risk bounds are
usually provable for estimators based on empirical model selection
criteria (see, e.g., \cite{BarronBirgeMassart}) or aggregation; see,
for example, \cite{RT12}. Specializing to the present situation, in
order to adapt over $\Ps_k$ as $k$ varies, one constructs LSEs over
each $\Ps_k$ and then either selects one estimator from this
collection by an empirical model selection criterion or aggregates
these estimators with data-dependent weights. In this particular
situation, such estimators are very difficult to compute as minimizing
the LS criterion over $\Ps_k$ is a nonconvex optimization problem. In
contrast, the LSE can be easily computed by a convex optimization
problem. It is remarkable that the LSE, which is constructed with no
explicit model selection criterion in mind, achieves adaptive risk
bound \eqref{hh}. In the next example, we illustrate this adaptation
for the LSE in convex regression.
\end{example}

\begin{example}[(Convex regression)]\label{conre}
Convex regression with equispaced design points corresponds to $\K_{r,
s}^n$ with $r = s = 1,
w_{-1} = w_1 = 1$ and $w_0 = -2$. It turns out that the statistical
dimension of this cone satisfies
%
\begin{equation}
\label{zb} \delta\bigl(\K_{1, 1}^n\bigr) \leq C \bigl(
\log(en)\bigr)^{5/4}\qquad \qt{for all $n \geq1$,}
\end{equation}
where $C$ is a universal positive constant. This is proved
in \cite{GSvex}, Theorem~3.1, via metric entropy results for classes
of convex functions.

Let $\alpha= (\alpha_1, \ldots, \alpha_n) \in\K_{-1, 1}^n$ with
$k(\alpha) = k$. Let $t_1, \ldots, t_k$ denote the values of $t$
where the inequality $2 \theta_t \leq\theta_{t-1} + \theta_{t+1}$
is strict. With $n_1 := t_1$, $n_i := t_i - t_{i-1}$ for $i = 2,
\ldots, k$ and $n_{k+1} = n - t_k$, we have, from \eqref{tde}
and \eqref{zb},
\[
\tau(\alpha) = \sum_{i=1}^{k+1} \delta\bigl(
\K_{1, 1}^{n_i}\bigr) \leq C \sum_{i=1}^{k+1}
(\log en_i )^{5/4}.
\]
Using the fact that $x \mapsto(\log x)^{5/4}$ is concave for $x
\geq e$, we have (note that $\sum_i n_i = n$)
\[
\tau(\alpha) \leq(k+1) \biggl(\log\frac{en}{k+1} \biggr)^{5/4}.
\]

Inequality \eqref{zens} then becomes
\[
\E_{\theta} \ell^2\bigl(\theta, \hat{\theta}\bigl(Y;
\K^{n}_{1, 1}\bigr)\bigr) \leq C \inf_{\alpha\in\K_{1,1}^n}
\biggl(\ell^2(\theta, \alpha) + \frac{\sigma^2 (k(\alpha) + 1)}{n} \biggl(\log
\frac{en}{k(\alpha) +
1} \biggr)^{5/4} \biggr).
\]
Note that the quantity $1 + k(\alpha)$ can be interpreted as the
number of affine pieces of the convex sequence $\alpha$. This risk
bound is the analogue of inequality \eqref{hh} for convex regression,
and it highlights the adaptation of the convex LSE to piecewise affine
convex functions. A weaker version of this inequality appeared in \cite{GSvex},
Theorem~2.3.
\end{example}

\subsection{Proof of Theorem \texorpdfstring{\protect\ref{zen}}{2.1}}
We now prove Theorem~\ref{zen}. We shall first state some general
results (Lemmas \ref{lili}, \ref{hol} and \ref{apps}), whose proofs
can be found in the supplementary material \cite{ChatEtAlSupp14}, for
the risk of $\E_{\theta} \ell^2(\theta, \hat{\theta}(Y;
\K))$, which hold for every $\K$ of the
form \eqref{bla}. Theorem~\ref{zen} will then be proved by
specializing these results for $\K= \K_{r, s}^n$.

We begin by recalling a result of Meyer and Woodroofe \cite{MW00} who
related the risk of $\hat{\theta}(Y; \K)$ to the function $D(\cdot;
\K)$. Specifically, \cite{MW00}, Proposition~2, proved that
\[
\E_0 \ell^2\bigl(0, \hat{\theta}(Y; \K)
\bigr) = \frac{\sigma^2\delta(\K
)}{n} = \frac{\sigma^2}{n} \E_0 D(Y; \K)
\]
and
%
\begin{equation}
\label{mw2} \E_{\theta} \ell^2\bigl(\theta, \hat{\theta}(Y;
\K)\bigr) \leq\frac
{\sigma^2}{n} \E_{\theta} D(Y; \K) \qquad\qt{for every $\theta
\in\K$}.
\end{equation}
These can be proved via Stein's lemma; see \cite{MW00}, proof of Proposition~2. It might be helpful to observe here that the function
$D(y; \K)$ satisfies $D(ty; \K) = D(y; \K)$ for every $t \in\R$, and
this is a consequence of the fact that $\hat{\theta}(ty; \K) =
t\hat{\theta}(y; \K)$ and the characterization \eqref{ctz}.

Our first lemma below states that the risk of the LSE is equal to
$\sigma^2\delta(\K)/n$ for all $\theta$ belonging to the lineality
space $\lin:= \{\theta\in\R^n \dvtx A \theta= 0\}$ of $\K$. The
lineality space $\lin$ will be crucial in the proof of
Theorem~\ref{zen}. The lineality space of the cone for isotonic
regression is the set of all constant sequences. The lineality space
of the cone for convex regression is the set of all affine sequences.
Also, we say that two convex polyhedral cones $\K_1$ and $\K_2$ are
orthogonal if $\langle\gamma_1,\gamma_2 \rangle= 0$ for all $\gamma
_1 \in\K_1$ and $\gamma_2 \in\K_2$.

\begin{lemma}\label{lili}
For every $\theta\in\R^n$ with $\theta= \gamma_1 + \gamma_2$ for
some $\gamma_1 \in\lin$ and $\gamma_2 \perp\K$ (i.e., $\langle
\gamma_2, \omega\rangle= 0$ for all $\omega\in\K$),
we have $\E_{\theta} D(Y; \K) = \delta(\K)$.
\end{lemma}


\begin{lemma}\label{hol}
Let $\K$ be an arbitrary convex polyhedral cone. Suppose $\K_1, \ldots
, \K_{l}$
are orthogonal polyhedral cones with lineality spaces $\lin_1, \ldots,
\lin_l$ such that $\K\subseteq\K_1 + \cdots+ \K_l$. Then
\[
\E_{\theta} D(Y; \K) \leq2\bigl(\delta(\K_1) + \cdots+
\delta(\K_l)\bigr)\qquad \qt{for every $\theta\in\K\cap(\lin_1
+ \cdots+ \lin_l)$}.
\]
\end{lemma}

The next lemma allows us to bound the risk of the LSE at $\theta$ by a
combination of the risk at $\alpha$ and the distance between $\theta$
and $\alpha$.

\begin{lemma}\label{apps}
The risk of the LSE satisfies the following inequality:
\[
\E_{\theta} \ell^2\bigl(\theta, \hat{\theta}(Y; \K)\bigr)
\leq3 \inf_{\alpha
\in\K} \bigl[ 2 \ell^2(\theta, \alpha) +
\E_{\alpha} \ell^2\bigl(\alpha, \hat{\theta}(Y; \K)\bigr) \bigr]\qquad
\qt{for every $\theta\in\K$}.
\]
\end{lemma}

We are now ready to prove Theorem~\ref{zen}.
\begin{pf*}{Proof of Theorem~\ref{zen}}
By Lemma~\ref{apps}, it is enough to prove that
\[
\E_{\alpha} \ell^2\bigl(\alpha, \hat{\theta}\bigl(Y;
\K_{r, s}^n\bigr)\bigr) \leq 2\bigl(1 + k(\alpha)\bigr)
\frac{\sigma^2\delta(\K_{r, s}^n)}{n} \qquad\qt{for every $\alpha\in\K_{r, s}^n$}.
\]
Fix $\alpha\in\K_{r, s}^n$, and let $k = k(\alpha)$, which means that
$k$ of the inequalities $\sum_{j=-r}^s w_j \alpha_{t+j} \geq0$ for $1
+ r \leq t \leq n-s$ are strict. Let $1+r \leq t_1 < \cdots< t_k \leq
n - s$ denote the indices of the inequalities that are strict. We
partition the set $\{1, \ldots, n\}$ into $k+1$ disjoint sets $E_0,
\ldots, E_{k}$ where
\[
E_0 := \{1, \ldots, t_1 - 1 + s\},\qquad  E_{k} :=
\{t_k + s, \ldots, n\}
\]
and
\[
E_i := \{t_i + s, \ldots, t_{i+1} - 1 + s\}\qquad
\qt{for $1 \leq i \leq k-1$}.
\]
Also for each $0 \leq i \leq k$, let
\[
F_i := \{t \in\mathbb{Z} \dvtx t-r \in E_i \mbox{ and }
t + s \in E_i \}.
\]
We now apply Lemma~\ref{hol} with
\[
\K_i := \Biggl\{\theta\in\R^n \dvtx
\theta_j = 0 \mbox{ for } j \notin E_i \mbox{ and } \sum
_{j=-r}^s w_j
\theta_{t + j} \geq0 \mbox{ for } t \in F_i \Biggr\}
\]
for $i = 0, \ldots, k$. The lineality space of $\K_i$ is, by
definition,
\[
\lin_i = \Biggl\{\theta\in\R^n \dvtx
\theta_j = 0 \mbox{ for } j \notin E_i \mbox{ and } \sum
_{j=-r}^s w_j
\theta_{t + j} = 0 \mbox{ for } t \in F_i \Biggr\}.
\]
$\K_0, \ldots, \K_k$ are orthogonal convex polyhedral cones because
$E_0, \ldots, E_k$ are disjoint. Also $\K\subseteq\K_0 + \cdots+
\K_k$ because every $\theta\in\K$ can be written as $\theta=
\sum_{i=0}^k \theta^{(i)}$ where $\theta_j^{(i)} := \theta_j I \{j
\in
E_i\}$ (it is easy to check that $\theta^{(i)} \in\K_i$ for each
$i$). Further, note that $\alpha\in\lin_0 + \cdots+ \lin_k$ since
$\alpha^{(i)} \in\lin_i$ for every $i$. Lemma~\ref{hol} thus gives
$\E_{\alpha} D(Y; \K) \leq2 \sum_{i=0}^k
\delta(\K_i)$. Inequality \eqref{mw2} then implies that
\[
\E_{\alpha} \ell^2\bigl(\alpha, \hat{\theta}(Y; \K)\bigr) \leq
\frac{2
\sigma^2}{n} \sum_{i=0}^k \delta(
\K_i).
\]
It is now easy to check that $\delta(\K_i) = \delta(\K_{r,
s}^{|E_i|})$ for each $i$ which proves \eqref{zens}. The proof
of \eqref{zen.eq} is now complete by the observation $\delta(\K_{r,
s}^{|E_i|}) \leq\delta(\K_{r, s}^n)$ as $|E_i| \leq n$.
\end{pf*}

\section{Risk bound in isotonic regression}\label{Isoreg}
In this section, we provide a proof of inequality \eqref{pyaa} using
an explicit formula of the LSE in isotonic regression. Our proof does
not require normality of $\varepsilon_1, \ldots, \varepsilon_n$. We also
explain the similarities between inequalities \eqref{zen.eq}
and \eqref{pyaa}. Before we get to inequality \eqref{pyaa}, however,
we give a brief overview of existing theoretical results in isotonic
regression.

Usually isotonic regression is posed as a function estimation problem
in which the unknown object of interest is a nondecreasing function
$f_0$, and one observes data from model \eqref{eq:Mdl2} with
$\theta_i := f_0(x_i), i = 1, \ldots, n$, where $x_1 < \cdots< x_n$ are
fixed design points. The most natural and commonly used estimator for
this problem is the monotone LSE defined as any
nondecreasing function $\hat{f}_{\mathrm{ls}}$ on $\R$ for which
$(\hat{f}_{\mathrm{ls}}(x_1), \ldots, \hat{f}_{\mathrm{ls}}(x_n)) = \hat{\theta}(Y;
\M)$. This estimator was proposed by \cite{Brunk55} and \cite
{AyerEtAl55}; also see \cite{G56} for the related problem of
estimating a nonincreasing density. Note that $\hat{\theta}(Y; \M)$
can be computed easily using the \textit{pool adjacent violators
algorithm}; see \cite{RWD88}, Chapter~1.

Existing theoretical results on isotonic regression can
be grouped into two categories: (1) results on the behavior of the LSE
at an interior point (which is sometimes known as local behavior), and
(2) results on the behavior of a global loss function measuring how
far $\hat{f}_{\mathrm{ls}}$ is from $f_0$.

Results on the local behavior are proved, among others,
in \cite{Brunk70,Wright81,G83,G85,CD99,Cator2011,H14}. Under
certain
regularity conditions on the unknown function $f_0$ near the interior point
$x_0$, it was proved in \cite{Brunk70} that $\hat{f}_{\mathrm{ls}}(x_0)$
converges to
$f_0(x_0)$ at the rate $n^{-1/3}$ and also characterized the limiting
distribution of $n^{1/3}(\hat{f}_{\mathrm{ls}}(x_0) - f_0(x_0))$. In the
related (nonincreasing) density estimation problem, the authors of
\cite{G85,CD99,H14} demonstrated that if the interior point $x_0$
lies on a flat stretch of the underlying function, then the LSE (which
is also the nonparametric maximum likelihood estimator, usually known
as the Grenander estimator) converges to a nondegenerate limit at
rate $n^{-1/2}$, and they characterized the limiting
distribution. In \cite{Cator2011}, Cator demonstrated that the rate of
convergence of
$\hat{f}_{\mathrm{ls}}(x_0)$ to $f_0(x_0)$ depends on the local behavior of
$f_0$ near $x_0$, and explicitly described this rate for each $f_0$. In
this sense, the LSE $\hat{f}_{\mathrm{ls}}$ adapts automatically to the
unknown function $f_0$. In \cite{Cator2011}, it was also proved that
the LSE is optimal for local behavior by establishing a local
asymptotic minimax lower bound.

Often in monotone regression, the interest is in the estimation of the
entire function $f_0$, as opposed to just its value at one fixed
point. In this sense, it is more appropriate to study the behavior of
$\hat{f}_{\mathrm{ls}}$ under a global loss function. The most natural and
commonly studied global loss function in this setup is
\[
L(f,g) := \frac{1}{n} \sum_{i=1}^n
\bigl(f(x_i) - g(x_i) \bigr)^2 =
\ell^2(\check{f}, \check{g}),
\]
where $\check{f} := (f(x_1), \ldots, f(x_n))$ and $\check{g} :=
(g(x_1), \ldots, g(x_n))$. Note that under this loss
function, the function estimation problem becomes exactly the same as
the sequence estimation problem described in \eqref{eq:Mdl2}, where
the goal is to estimate the vector $\theta:= (\theta_1, \ldots,
\theta_n)$ under the constraint $\theta\in\M$ and the loss
function \eqref{eq:loss}. The behavior of $\hat{\theta}(Y; \M)$,
under the loss $\ell^2$, has been studied in a number of papers
including \cite{vdG90,vdG93,Donoho91,BM93,Wang96,MW00,Zhang02}. If one
looks at the
related (nonincreasing) density estimation problem, Birg{\'e} \cite{Birge89}
developed nonasymptotic risk bounds for the Grenander estimator,
measured with the $L_1$-loss, whereas Van de Geer \cite{vdG93} has
results on the
Hellinger distance. As mentioned in the \hyperref[sec1]{Introduction}, the strongest
existing bound on $\E_{\theta}\ell^2(\theta, \hat{\theta}(Y; \M
))$ is
due to \cite{Zhang02}, Theorem~2.2. We recalled this inequality
in \eqref{motw} and compared it to our bound \eqref{pyaa} in some
situations.

In the following theorem, we prove inequality \eqref{pyaa} using
explicit characterization of $\hat{\theta}(Y; \M)$ without requiring
normality of $\varepsilon_1, \ldots, \varepsilon_n$. In fact, we prove an
inequality that is slightly stronger than \eqref{pyaa}.

We need the following notation. For simplicity, we use $\hat{\theta}$
for $\hat{\theta}(Y; \M)$ and $\hat{\theta}_j$ for the components of
$\hat{\theta}(Y; \M)$. For any sequence $(a_1,a_2,\ldots,a_n) \in\R^n$
and any $1 \leq u \leq v \leq n$, let
%
\begin{equation}
\label{eq:ShortHand} \bar{a}_{u, v} := \frac{1}{v-u+1} \sum
_{j=u}^v a_j.
\end{equation}
We will use this notation mainly when $a$ equals $Y$, $\theta$ or
$\varepsilon$. Our proof uses ideas similar to those in \cite{Zhang02},
Section~2 (see Remark~\ref{annr} for details about the connections
with \cite{Zhang02}, Section~2) and is based on the following
explicit representation of the LSE $\hat{\theta}$ (see \cite{RWD88},
Chapter~1):
%
\begin{equation}
\label{eq:Charac} \hat{\theta}_j = \min_{v \geq j} \max
_{u \leq j} \bar{Y}_{u,v}.
\end{equation}
For $x \in\R$, we write $x_+ := \max
\{0,x\}$ and $x_- := - \min\{0,x\}$.
For $\theta\in\M$ and $\pi= (n_1, \ldots, n_{k+1}) \in\Pi$, let
\[
D_{\pi} (\theta) = \Biggl( \frac{1}{n} \sum
_{i =
1}^{k+1} \sum_{j=s_{i-1} + 1}^{s_i}
( \theta_j - \bar{\theta}_{s_{i-1} + 1,s_i} )^2
\Biggr)^{1/2},
\]
where $s_0 = 0$ and $s_i = n_1 + \cdots+ n_i$, for $1 \leq i \leq
k+1$. Like $V_{\pi}(\theta)$, this quantity $D_{\pi}(\theta)$ can also
be treated as a measure of the variation of $\theta$ with respect to
$\pi$. This measure also satisfies $D_{\pi_{\theta}}(\theta) = 0$ for
every $\theta\in\M$. Moreover,
\[
D_{\pi}(\theta) \leq V_{\pi}(\theta)\qquad \qt{for every $\theta\in
\M$ and $\pi\in\Pi$}.
\]
When $\pi= (n)$ is the trivial partition, $D_{\pi}(\theta)$ turns
out to be just the standard deviation of $\theta$. In general,
$D_{\pi}^2(\theta)$ is analogous to the within group sum of squares
term in ANOVA with the blocks of $\pi$ being the groups. Below, we
prove a stronger version of \eqref{pyaa} with $D_{\pi}(\theta)$
replacing $V_{\pi}(\theta)$ in the definition of $R(n; \theta)$.

\begin{theorem}\label{rafi}
Suppose $Y_1, \ldots, Y_n$ are observations from model \eqref{eq:Mdl2}
with $\varepsilon_1, \ldots, \varepsilon_n$ being i.i.d. with mean zero and
variance $\sigma^2$. For every $\theta\in\M$, the
risk of $\hat{\theta} = \hat{\theta}(Y; \M)$ satisfies the following
inequality:
%
\begin{equation}
\label{rafi.eq} \E_{\theta} \ell^2(\theta, \hat{\theta}) \leq4
\inf_{\pi\in
{\Pi}} \biggl(D^2_{\pi} (\theta) +
\frac{4 \sigma^2 (1 + k(\pi))}{n} \log \frac{en}{1 + k(\pi)} \biggr).
\end{equation}
\end{theorem}

\begin{pf}
Fix $1 \leq j \leq n$ and $0 \leq m \leq n-j$. By \eqref{eq:Charac},
we have
\[
\hat{\theta}_j = \min_{v \geq j} \max
_{u \leq j} \bar{Y}_{u, v} \leq\max_{u \leq j}
\bar{Y}_{u, j+m} = \max_{u \leq j} ( \bar{
\theta}_{u, j+m} + \bar{\varepsilon}_{u, j+m} ),
\]
where, in the last equality, we used $\bar{Y}_{u, v} = \bar{\theta}_{u,
v} + \bar{\varepsilon}_{u, v}$. By the monotonicity of $\theta$, we
have $\bar{\theta}_{u, j+m} \leq\bar{\theta}_{j, j+m}$ for all $u
\leq j$. Therefore, for every $\theta\in\M$, we get
\[
\hat{\theta}_j - \theta_j \leq(\bar{
\theta}_{j, j+m} - \theta_j) + \max_{u \leq j}
\bar{\varepsilon}_{u, j+m}.
\]
Taking positive parts, we have
\[
( \hat{\theta}_j - \theta_j )_+ \leq(\bar{\theta
}_{j, j+m} - \theta_j) + \max_{u \leq j} (\bar{
\varepsilon}_{u, j+m} )_+.
\]
Squaring and taking expectations on both sides, we obtain
\[
\E_{\theta} (\hat{\theta}_j - \theta_j
)_+^2 \leq\E _{\theta} \Bigl((\bar{\theta}_{j, j+m} -
\theta_j) + \max_{u \leq j} ( \bar{\varepsilon}_{u, j+m}
)_+ \Bigr)^2.
\]
Using the elementary inequality $(a+b)^2 \leq2a^2 + 2b^2$ we get
\[
\E_{\theta} (\hat{\theta}_{j} - \theta_j
)_+^2 \leq 2 (\bar{\theta}_{j, j+m} - \theta_j
)^2 + 2 {\E\max_{u
\leq j} ( \bar{\varepsilon}_{u, j+m}
)_+^2 }.
\]
We observe now that, for fixed integers $j$ and $m$, the process
$\{\bar{\varepsilon}_{u, j+m}, u = 1, \ldots, j\}$ is a martingale with
respect to the filtration $\F_1, \ldots, \F_j$ where $\F_i$ is the
sigma-field generated by the random variables $\varepsilon_1, \ldots,
\varepsilon_{i-1}$ and $\bar{\varepsilon}_{i, j+m}$. Therefore, by Doob's
inequality for submartingales (see, e.g., Theorem~5.4.3 of
\cite{Durrett10}), we have
\[
\E\max_{u \leq j} (\bar{\varepsilon}_{u, j+m}
)_+^2 \leq4 \E (\bar{\varepsilon}_{j, j+m} )_+^2
\leq4 \E (\bar{\varepsilon}_{j, j+m} )^2 \leq
\frac{4 \sigma^2}{m + 1}.
\]
So using the above result we get the following pointwise upper bound
for the
positive part of the risk:
%
\begin{equation}
\label{eq:sc1} \E_{\theta} (\hat{\theta}_{j} -
\theta_j )_+^2 \leq 2 (\bar{\theta}_{j, j+m} -
\theta_j )^2 + \frac{8
\sigma^2}{m + 1}.
\end{equation}
Note that the above upper bound holds for any arbitrary $m$, $0 \leq m
\leq n-j$.
By a similar argument we can get the following pointwise upper bound
for the negative part of risk which now holds for any $m$, $0 \leq m
\leq j$:
%
\begin{equation}
\label{eq:sc2} \E_{\theta} (\hat{\theta}_{j} -
\theta_j )_-^2 \leq 2 (\theta_j - \bar{
\theta}_{j-m, j} )^2 + \frac{8 \sigma
^2}{m + 1}.
\end{equation}
Let us now fix $\pi= (n_1, \ldots, n_{k+1}) \in\Pi$. Let $s_0 := 0$
and $s_i := n_1 + \cdots+ n_i$ for $1 \leq i \leq k+1$. For each $j
=1,\ldots, n$, we define two integers $m_1(j)$ and $m_2(j)$ in the
following way: $m_1(j) = s_i - j$ and $m_2(j) = j - 1 - s_{i-1}$ when
$s_{i-1} + 1 \leq j \leq s_i$. We use this choice of $ m_1(j) $
in \eqref{eq:sc1} and $ m_2(j) $ in \eqref{eq:sc2} to obtain
$\E_{\theta} (\hat{\theta}_j - \theta_j)^2 \leq A_j +
B_j$ where
\[
A_j := 2 (\bar{\theta}_{j, j+m_1(j)} - \theta_j
)^2 + \frac{8 \sigma^2}{m_1(j) + 1}
\]
and
\[
B_j := 2 (\theta_j - \bar{\theta}_{j-m_2(j), j}
)^2 + \frac{8 \sigma^2}{m_2(j) + 1}.
\]
This results in the risk bound
\[
\E_{\theta} \ell^2(\theta,\hat{\theta}) \leq
\frac{1}{n} \sum_{j=1}^n
A_j + \frac{1}{n} \sum_{j=1}^n
B_j.
\]
We shall now prove that
%
\begin{equation}
\label{kk1} \frac{1}{n} \sum_{j=1}^n
A_j \leq2 D^2_{\pi}(\theta) +
\frac{8(k+1)\sigma^2}{n} \log\frac{en}{k+1}
\end{equation}
and
%
\begin{equation}
\label{kk2} \frac{1}{n} \sum_{j=1}^n
B_j \leq2 D^2_{\pi}(\theta) +
\frac{8(k+1)\sigma^2}{n} \log\frac{en}{k+1}.
\end{equation}
We give below the proof of \eqref{kk1}, and the proof of \eqref{kk2} is
nearly identical. Using the form of $A_j$, we break up $\frac{1}{n}
\sum_{j=1}^n A_j$ into two terms. For the first term, note that $j +
m_1(j) = s_i$, for $s_{i-1}+1 \leq j \leq s_i$ and therefore
\[
\sum_{j=1}^n (\bar{\theta}_{j, j+m_1(j)}
- \theta_j )^2 = \sum_{i=1}^{k+1}
\sum_{j=s_{i-1} + 1}^{s_i} (\bar{\theta}_{j,
s_i}
- \theta_j )^2.
\]
By Lemma~11.2 in the supplementary material \cite
{ChatEtAlSupp14}, we get
\[
\sum_{j=s_{i-1} + 1}^{s_i} (\bar{\theta}_{j, s_i}
- \theta_j )^2 \leq\sum_{j=s_{i-1} + 1}^{s_i}
(\bar{\theta}_{s_{i-1}+1,
s_i} - \theta_j )^2
\]
for every $i = 1, \ldots, k+1$. Thus summing over $i = 1, \ldots, k+1$,
and multiplying by $2/n$ proves that the first term in $\frac{1}{n}
\sum_{j=1}^n A_j$ is bounded from above by $2 D^2_{\pi}(\theta)$. To
bound the second term, we write
%
\begin{equation}
\label{eq:2ndTerm} \sum_{j=1}^n
\frac{1}{m_1(j) + 1} = \sum_{i=1}^{k+1} \sum
_{j =
s_{i-1}+1}^{s_i} \frac{1}{s_i - j + 1} = \sum
_{i=1}^{k+1} \biggl(1 + \frac{1}{2} +
\cdots+ \frac{1}{n_i} \biggr).
\end{equation}
Since the harmonic series $\sum_{i=1}^l 1/l$ is at most $1 + \log l$
for $l \geq1$, we obtain
\[
\sum_{j=1}^n \frac{1}{m_1(j)+1} \leq k+1
+ \sum_{i=1}^{k+1} \log n_i \leq
k+1 + (k+1) \log \biggl(\frac{\sum_i n_i}{k+1} \biggr),
\]
where the last inequality is a consequence of the concavity of the
logarithm function. This proves \eqref{kk1} because $\sum_{i} n_i =
n$. Combining \eqref{kk1} and \eqref{kk2} proves the theorem.
\end{pf}

\begin{remark}\label{remp}
For each $\pi= (n_1, \ldots, n_{k+1}) \in\Pi$, let $\M_{\pi}$ denote
the set of all $\alpha\in\M$ such that $\alpha$ is constant on
each set $\{j \dvtx s_{i-1} + 1 \leq j \leq s_i\}$ for $i = 1, \ldots,
k+1$. Then it is easy to see that
\[
\inf_{\alpha\in\M_{\pi}} \ell^2(\theta, \alpha) =
D_{\pi}^2(\theta).
\]
Using this, it is easy to see that inequality \eqref{rafi.eq} is
equivalent to
%
\begin{equation}
\label{nsv} \E_{\theta} \ell^2(\theta, \hat{\theta}) \leq4
\inf_{\alpha\in
\M} \biggl(\ell^2(\theta, \alpha) +
\frac{4 \sigma^2 (1 +
k(\alpha))}{n} \log\frac{en}{1 + k(\alpha)} \biggr).
\end{equation}
Inequality \eqref{rafi.eq} therefore differs from
inequality \eqref{hh} only by its multiplicative constants. It should be
noted that we proved \eqref{hh} assuming normality of $\varepsilon_1,
\ldots, \varepsilon_n$ while \eqref{rafi.eq} was proved without
using normality. Inequality \eqref{pyaa} is slightly weaker
than \eqref{rafi.eq} because $D_{\pi}(\theta) \leq
V_{\pi}(\theta)$. We still work with \eqref{pyaa} in isotonic
regression as opposed to \eqref{rafi.eq} because it is easier to
compare \eqref{pyaa} to existing inequalities, and also, as we shall
show in Section~\ref{lmnop}, inequality \eqref{pyaa} is nearly
optimal.
\end{remark}

\begin{remark}\label{annr}
Bounding the infimum in the right-hand side of \eqref{rafi.eq} by
taking $\pi= \pi_{\theta}$ and letting $k(\theta) = k$, we obtain
%
\begin{equation}
\label{annr.eq} \E_{\theta} \ell^2(\theta, \hat{\theta}) \leq
\frac{16 \sigma^2
(1 + k)}{n} \log\frac{en}{1 + k}.
\end{equation}
This inequality might be implicit in the arguments of \cite{Zhang02},
Section~2. It might be possible to prove \eqref{annr.eq} by
applying \cite{Zhang02}, Theorem~2.1, to each of the $k+1$ constant
pieces of $\theta$ and by bounding the resulting quantities via
arguments in \cite{Zhang02}, proof of Theorem~2.2.
\end{remark}

\begin{remark}\label{ach}
Note that $k(\theta)$ does not have to be small for \eqref{pyaa} to be
an improvement of \eqref{motw}. One only needs that $V_{\pi}(\theta)$
be small for some partition $\pi$ with small $k(\pi)$. Equivalently,
from \eqref{nsv}, one needs that $\ell^2(\theta, \alpha)$ is small for
some $\alpha\in\M$ with small $k(\alpha)$. This is illustrated below.

Let $\{a_j, j \geq1\}$ be an arbitrary countable subset of $[0, 1]$,
and let $\{p_j, j \geq1\}$ denote any probability sequence, that is, $p_j
\geq0, \sum_j p_j = 1$. Fix $n \geq1$, and let $\theta_i := \sum_{j\dvtx a_j \leq i/n} p_j$ for $i = 1, \ldots, n$. We will argue below that,
for many choices of $\{p_j, j \geq1\}$, inequality \eqref{nsv} gives
a faster rate of convergence than $n^{-2/3}$, even though $k(\theta)$
can be as large as $n$.

Indeed, fix $1 \leq k \leq n$, and define $\alpha_i = \sum_{j \leq k\dvtx a_j \leq i/n} p_j$. It is then clear that $k(\alpha) \leq k$. Also
for each $1 \leq i \leq n$, we have
\[
0 \leq\theta_i - \alpha_i = \sum
_{j > k\dvtx a_j \leq i/n} p_j \leq \sum
_{j > k} p_j.
\]
This implies that $\ell(\theta, \alpha) \leq\sum_{j\dvtx j > k}
p_j$. Thus inequality \eqref{nsv} gives
%
\begin{equation}
\label{alit} \E_{\theta} \ell^2(\theta, \hat{\theta}) \leq4
\inf_{k} \biggl[ \biggl( \sum_{j\dvtx j > k}
p_j \biggr)^2 + \frac{4 \sigma^2(1 + k)}{n} \log
\frac{en}{1 + k} \biggr].
\end{equation}
When $\sum_{j\dvtx j > k} p_j = o(k^{-1})$, it can be checked that the
bound above is faster than $n^{-2/3}$. This happens, for
instance, when $p_j \propto j^{-a}$ for $a \geq3$. In fact, when $p_j =
2^{-j}$, \eqref{alit} gives the parametric rate up to logarithmic
factors.

However, when, say $p_j \propto j^{-a}$ for $1 < a \leq
2$, \eqref{alit} does not give a rate that is faster than
$n^{-2/3}$. It might be possible here to use a different approximation
vector $\alpha\in\M$ which would still yield a rate of
$o(n^{-2/3})$, but we
do not have a proof of this. A result from the literature that is
relevant here is
\cite{Zhang02}, inequality (2.10). For the vectors $\theta$
considered above (and in certain more general situations), this
inequality gives an asymptotic bound of $o(n^{-2/3})$ (without
quantifying the exact order) for the risk for all choices of $\{a_j\}$
and $\{p_j\}$, even for $p_j \propto j^{-a}, 1 < a \leq2$.
\end{remark}
%
\begin{example}
We prove in Theorem~\ref{sabya2} in the next section that the bound
given by Theorem~\ref{rafi} is always smaller than a logarithmic
multiplicative factor of the usual cube root rate of convergence for
every $\theta\in\M$ with $V(\theta) > 0$. Here, we shall
demonstrate this in the special case of the sequence $\theta= (1/n, 2/n,
\ldots, 1)$ where the bound in \eqref{rafi.eq} can be calculated
exactly. Indeed, if $\pi= (n_1, \ldots, n_k)$ with $n_i \geq1$ and
$\sum_{i=1}^k n_i = n$, direct calculation gives
\[
D^2_{\pi}(\theta) = \frac{1}{12n^3} \Biggl(\sum
_{i=1}^k n_i^3 - n \Biggr).
\]
Now H\"{o}lder's inequality gives $n = \sum_{i=1}^k n_i \leq
(\sum_{i=1}^k n_i^3)^{1/3} k^{2/3}$ which means that $\sum_{i=1}^k
n_i^3 \geq n^3/k^2$. Therefore, for every fixed $k \in\{1, \ldots,
n\}$ such that $n/k$ is an integer, $D_{\pi}^2(\theta)$ is minimized
over all partitions $\pi$ with $k(\pi) = k$ when $n_1 = n_2 = \cdots=
n_k = n/k$. This gives $\inf_{\pi\dvtx  k(\pi) = k} D^2(\pi) =
(k^{-2} - n^{-2})/12$. As a consequence, Theorem~\ref{rafi} yields the
bound
\[
\E_{\theta} \ell^2(\theta, \hat{\theta}) \leq
\frac{1}{3} \inf_{k\dvtx n/k \in\mathbb{Z}} \biggl(\frac{1}{k^2} -
\frac{1}{n^2} + \frac{48
\sigma^2 k}{n} \log(en/k) \biggr).
\]
Now with the choice $k \sim(n/\sigma^2)^{1/3}$, we get the cube root
rate for $\hat{\theta}$ up to logarithmic multiplicative factors in
$n$. We generalize this to arbitrary $\theta\in\M$ with $V(\theta) >
0$ in Theorem~\ref{sabya2}.
\end{example}
%
\section{The quantity \texorpdfstring{$R(n;\theta)$}{R(n;theta)}}\label{hd}
In this section, we state some results about the quantity $R(n;
\theta)$ appearing in our risk bound \eqref{pyaa}. Recall also the
quantity $R_Z(n \theta)$ that appears in \eqref{motw}. The first
result of this section states that $R(n;\theta)$ is always bounded
from above by $R_Z(n; \theta)$ up to a logarithmic multiplicative
factor in $n$. This implies that \eqref{pyaa} is always only slightly
worse off
than \eqref{motw} (by a logarithmic multiplicative factor) while being
much better when $\theta$ is well-approximated by some $\alpha\in\M$
for which $k(\alpha)$ is small. Recall that $V(\theta) :=
\theta_n - \theta_1$. The proofs of all the results in this section
can be found in the supplementary material \cite{ChatEtAlSupp14}.

\begin{theorem}\label{sabya2}
For every $\theta\in\M$, we have
%
\begin{equation}
\label{sabya2.eq} R(n;\theta) \leq16 \log(4n) \biggl(\frac{\sigma^2
V(\theta)}{n}
\biggr)^{2/3}
\end{equation}
whenever
%
\begin{equation}
\label{sabya2.con} n \geq\max \biggl(2, \frac{8 \sigma^2}{V^2(\theta)}, \frac{V(\theta)}{\sigma}
\biggr).
\end{equation}
\end{theorem}

In the next result, we characterize $R(n;\theta)$ for certain strictly
increasing sequences $\theta$ where we show that it is essentially of the
order $(\sigma^2 V(\theta)/n)^{2/3}$. In some sense, $R(n; \theta)$ is
maximized for these strictly increasing sequences. The prototypical
sequence we have in mind here is $\theta_i = i/n$ for $1 \leq i \leq
n$.

\begin{theorem}\label{cst}
Suppose $\theta_1 < \theta_2 < \cdots<\theta_n$ with
%
\begin{equation}
\label{lcr} \min_{2 \leq i \leq n} (\theta_i -
\theta_{i-1} ) \geq \frac{c_1 V(\theta)}{n}
\end{equation}
for a positive constant $c_1 \leq1$. Then we have
%
\begin{equation}
\label{cst.eq} 12 \biggl( \frac{c_1\sigma^2 V(\theta)}{n} \biggr)^{2/3} \leq R(n,
\theta) \leq16 \biggl( \frac{\sigma^2 V(\theta)}{n} \biggr)^{2/3} \log(4n)
\end{equation}
provided
%
\begin{equation}
\label{cst.con} n \geq\max \biggl(2, \frac{8 \sigma^2}{V^2(\theta)}, \frac{2V(\theta)}{\sigma}
\biggr).
\end{equation}
\end{theorem}
%
\begin{remark}
An important situation where \eqref{lcr} is satisfied is when $\theta$
arises from sampling a function on $[0, 1]$ at the points $i/n$ for
$i = 1, \ldots, n$, assuming that the derivative of the function is
bounded from below by a positive constant.
\end{remark}
Next we describe sequences $\theta$ for which $R(n;\theta)$ is
$(k(\theta)\sigma^2/n) \log(en/k(\theta))$, up to multiplicative
factors. For these sequences our risk bound is potentially far
superior to $R_Z(n;\theta)$.

\begin{theorem}\label{vt}
Let $k = k(\theta)$ with $\{y\dvtx y = \theta_j \mbox{ for some } j \} =
\{\theta_{0,1}, \ldots, \theta_{0,k} \}$ where $\theta_{0,1} <
\cdots
< \theta_{0,k}$. Then
%
\begin{equation}
\label{vt.eq} \frac{\sigma^2 k}{n} \log\frac{en}{k} \leq R(n;\theta) \leq
\frac{16
\sigma^2 k}{n} \log\frac{en}{k}
\end{equation}
provided
%
\begin{equation}
\label{vt.con} \min_{2 \leq i \leq k} (\theta_{0,i} -
\theta_{0,i-1} ) \geq \sqrt{\frac{k\sigma^2}{n} \log\frac{en}{k}}.
\end{equation}
\end{theorem}

\section{Local minimax optimality of the LSE}\label{lmnop}
In this section, we establish an optimality property of the
LSE. Specifically, we show that $\hat{\theta}$ is locally minimax
optimal in a nonasymptotic sense. ``Local'' here refers to a ball
$\{t\dvtx \ell_{\infty}^2(t, \theta) \leq c R(n; \theta) \}$ around the
true parameter $\theta$ for a positive constant $c$. The reason we
focus on local minimaxity, as opposed to the more traditional notion of
global minimaxity, is that the rate $R(n; \theta)$ changes with
$\theta$. Note that, moreover, lower bounds on the global minimax risk
follow from our local minimax lower bounds. Such an optimality theory
based on local minimaxity has been pioneered by Cai and Low \cite{CaiLowFwork}
and Cator \cite{Cator2011} for the problem of estimating a convex or
monotone function at a point.

We start by proving an upper bound for the local supremum
risk of $\hat{\theta}$. Recall that $\ell_{\infty}(t, \theta) :=
\max_{1 \leq i \leq n} |t_i - \theta_i|$.

\begin{lemma}\label{lub.eq}
The following inequality holds for every $\theta\in\M$ and $c > 0$:
%
\begin{equation}
\label{eq:lub.eq} \sup_{t \in\M\dvtx \ell^2_{\infty}(t, \theta) \leq c R(n; \theta)} \E_t
\ell^2(t, \hat{\theta}) \leq2(1+4c) R(n;\theta).
\end{equation}
\end{lemma}

\begin{pf}
Inequality \eqref{pyaa} gives $\E_t \ell^2(t, \hat{\theta}) \leq
R(n;t)$ for every $t \in\M$. Fix $\pi\in\Pi$. By the triangle
inequality, we get $V_{\pi}(t) \leq2 \ell_{\infty}(t, \theta) +
V_{\pi}(\theta)$. As a result, whenever $\ell_{\infty}(t, \theta)
\leq c R(n; \theta)$, we obtain
\[
V^2_{\pi}(t) \leq2 V_{\pi}^2(\theta) +
8 \ell_{\infty}^2(t, \theta) \leq2 V_{\pi}^2(
\theta) + 8 c R(n; \theta).
\]
As a consequence,
\begin{eqnarray*}
\E_t \ell^2(t, \hat{\theta}) \leq R(n; t) &\leq& \inf
_{\pi\in\Pi} \biggl(2 V_{\pi}^2(\theta) +
\frac{4 \sigma^2 k(\pi)}{n} \log \frac{n}{k(\pi)} \biggr) + 8c R(n; \theta)
\\
&\leq&2 R(n; \theta) + 8 c R(n; \theta).
\end{eqnarray*}
This proves \eqref{eq:lub.eq}.
\end{pf}
We now show that $R(n; \theta)$, up to logarithmic factors in
$n$, is a lower bound for the local minimax risk at $\theta$, defined as
the infimum of the right-hand side of~\eqref{eq:lub.eq} over all
possible estimators $\hat{\theta}$. We prove this under each of the
assumptions (1) and (2) (stated in the \hyperref[sec1]{Introduction}) on
$\theta$. Specifically, we prove the two inequalities \eqref{inmi1}
and \eqref{inmi2}. These results mean
that, when $\theta$ satisfies either of the two assumptions (1)
or (2), no estimator can have a
supremum risk significantly better than $R(n;\theta)$ in the local
neighborhood $\{t \in\M\dvtx \ell_{\infty}^2(t,\theta) \lesssim R(n;
\theta)\}$. On the other hand, Lemma~\ref{lub.eq} states that the
supremum risk of the LSE over the same local neighborhood is bounded
from above by a constant multiple of $R(n;\theta)$. Putting these two
results together, we deduce that the LSE is approximately locally
nonasymptotically minimax for such sequences $\theta$. We use the
qualifier ``approximately'' here because of the presence of logarithmic
factors on the right-hand sides of \eqref{inmi1} and \eqref{inmi2}.

We make here the assumption that the errors $\varepsilon_1, \ldots,
\varepsilon_n$ are independent and normally distributed with mean zero
and variance $\sigma^2$. For each $\theta\in\M$, let $\P_\theta$
denote the
joint distribution of the data $Y_1, \ldots, Y_n$ when the true
sequence equals~$\theta$. As a consequence of the normality of the
errors, we have
\[
D(\P_\theta\|\P_t) = \frac{n}{2 \sigma^2}
\ell^2(t, \theta),
\]
where $D(P\|Q)$ denotes the Kullback--Leibler divergence between the
probability measures $P$ and $Q$. Our main tool for the proofs is
Assouad's lemma, the following version of which is a consequence of
Lemma~24.3 of \cite{vaart98book}, page 347.

\begin{lemma}[(Assouad)]\label{suad}
Let $m$ be a positive integer and suppose that, for each $\tau\in
\{0, 1\}^m$, there is an associated nondecreasing sequence $\theta
^{\tau}$ in $N(\theta)$, where $N(\theta)$ is a neighborhood of
$\theta$. Then the following inequality holds:
\[
\inf_{\hat{t}} \sup_{t \in N(\theta)} \E_t
\ell^2(t,\hat{t}) \geq\frac{m}{8} \min_{\tau\neq\tau'}
\frac{\ell^2(\theta^{\tau}, \theta^{\tau'})}{\ham(\tau, \tau')}
 \min_{\ham(\tau, \tau') = 1}\bigl (1 - \|\P_{\theta^{\tau}} -
\P_{\theta^{\tau'}}\|_{\mathrm{TV}} \bigr),
\]
where $\ham(\tau, \tau') := \sum_{i} I\{\tau_i \neq\tau'_i\}$ denotes
the Hamming distance between $\tau$ and $\tau'$ and $\|\cdot\|_{\mathrm{TV}}$
denotes the total variation distance between probability measures. The
infimum here is over all possible estimators $\hat{t}$.
\end{lemma}

Inequalities \eqref{inmi1} and \eqref{inmi2} are proved in the
next two subsections.

\subsection{Uniform increments} \label{ui}
In this section, we assume that $\theta$ is a strictly increasing
sequence with $V(\theta) = \theta_n - \theta_1 > 0$ and that
%
\begin{equation}
\label{cr} \frac{c_1 V(\theta) }{n} \leq\theta_i -
\theta_{i-1} \leq\frac{c_2
V(\theta)}{n} \qquad\qt{for $i = 2, \ldots, n$}
\end{equation}
for some $c_1 \in(0, 1]$ and $c_2 \geq1$. Because $V(\theta) =
\sum_{i=2}^n (\theta_i - \theta_{i-1} )$, assumption \eqref{cr} means
that the increments of $\theta$ are in a sense uniform. An important
example in which \eqref{cr} is satisfied is when $\theta_i = f_0(i/n)$
for some function $f_0$ on $[0, 1]$ whose derivative is uniformly
bounded from above and below by positive constants.

In the next theorem, we prove that the local minimax risk at $\theta$
is bounded from below by $R(n;\theta)$ (up to logarithmic
multiplicative factors) when $\theta$ satisfies \eqref{cr}.

\begin{theorem}\label{oval}
Suppose $\theta$ satisfies \eqref{cr}, and let
\[
\nbr(\theta) := \biggl\{t \in\M\dvtx \ell^2_{\infty}(t,
\theta) \leq \biggl(\frac{3c_2}{c_1} \biggr)^{2/3} \frac{R(n;\theta)}{12}
\biggr\}.
\]
Then the local minimax risk $\mix_n(\theta) := \inf_{\hat{t}} \sup_{t
\in\nbr(\theta)} \E_t \ell^2(t, \hat{t})$ satisfies the following
inequality:
%
\begin{equation}
\label{oval.eq} \mix_n(\theta) \geq \frac{c_1^2 3^{2/3}}{256
c_2^{4/3}} \biggl(
\frac{\sigma^2 V(\theta)}{n} \biggr)^{2/3} \geq \frac{c_1^2 3^{2/3}}{4096 c_2^{4/3}}
\frac{R(n;\theta)}{\log(4n)},
\end{equation}
provided
%
\begin{equation}
\label{oval.con} n \geq\max \biggl(2, \frac{24 \sigma^2}{V^2(\theta)}, \frac{2c_2
V(\theta)}{\sigma}
\biggr).
\end{equation}
\end{theorem}

Theorem~\ref{oval} is closely connected to minimax lower bounds for
Lipschitz classes of functions. Indeed, using the notation $v_f :=
(f(1/n), \ldots, f(1))$ for functions $f$ on $[0, 1]$, it can be argued
that
\[
\bigl\{t = \theta+ v_f \dvtx \bigl\|f'\bigr\|_{\infty}
\leq c_1 V(\theta), \|f\|_{\infty} \leq c_1'n^{-1/3}
\bigr\}
\]
is a subset of $\nbr(\theta)$ for appropriate positive constants
$c_1$ and $c_1'$. Lower bound \eqref{oval.eq} then follows from
usual lower bounds for Lipschitz classes which are
outlined, for example, in \cite{Tsybakovbook}, Chapter~2. A direct
proof of Theorem~\ref{oval} is included in the supplementary
material \cite{ChatEtAlSupp14}.

\subsection{Piecewise constant}\label{pc}
Here, we again show that the local minimax risk at $\theta$ is bounded
from below by $R(n; \theta)$ (up to logarithmic multiplicative
factors). The difference from the previous section is that we work
under a different assumption from \eqref{cr}. Specifically, we assume
that $k(\theta) = k$ and that the $k$ values of $\theta$ are
sufficiently well separated and prove inequality \eqref{inmi2}.

Let $k = k(\theta)$. There exist integers $n_1, \ldots, n_k$ with $n_i
\geq1$ and $n_1 + \cdots+ n_k = n$ such that $\theta$ is constant on
each set $\{j\dvtx s_{i-1}+1 \leq j \leq s_i \}$ for $i = 1, \ldots, k$
where $s_0 := 0$ and $s_i := n_1 + \cdots+ n_i$. Also, let the values
of $\theta$ on the sets $\{j\dvtx s_{i-1}+ 1 \leq
j \leq s_i \}$ for $i=1, \ldots, k$ be denoted by $\theta_{0,1} <
\cdots< \theta_{0,k}$.

\begin{theorem}\label{fani}
Suppose $c_1 n/k \leq n_i \leq c_2 n/k$ for all $1 \leq i \leq
k$ for some $c_1 \in(0, 1]$ and $c_2 \geq1$ and that
%
\begin{equation}
\label{masa} \min_{2 \leq i \leq k} (\theta_{0, i} -
\theta_{0, i-1} ) \geq \sqrt{\frac{k\sigma^2}{n} \log\frac{en}{k}}.
\end{equation}
Then, with $\nbr(\theta)$ defined as $ \{t \in\M\dvtx \ell_{\infty}^2(t,\theta) \leq R(n; \theta)  \}$, the local
minimax risk, $\mix_n(\theta) = \inf_{\hat{t}} \sup_{t \in
\nbr(\theta)} \E_t \ell^2(t, \hat{t})$, satisfies
\[
\mix_n(\theta) \geq\frac{c_1^{7/3}}{2^{31/3} c_2^2} R(n;\theta) \biggl(
\log\frac{en}{k} \biggr)^{-2/3},
\]
provided
%
\begin{equation}
\label{tie.con} \frac{n}{k} \geq\max \biggl( \biggl(\frac{4}{c_1^2} \log
\frac{en}{k} \biggr)^{1/3}, \exp \biggl(\frac{1-4c_1}{4c_1} \biggr)
\biggr).
\end{equation}
\end{theorem}
\begin{pf}
For notational convenience, we write
\[
\beta_n^2 := \frac{k\sigma^2}{n} \log\frac{en}{k}.
\]
First note that under assumption \eqref{masa}, Theorem~\ref{vt}
implies that $\beta_n^2 \leq R(n; \theta)$.

Let $1 \leq l \leq\min_{1 \leq i \leq k} n_i$ be a positive integer
whose value will be specified later, and let $m_i := \lfloor n_i/l
\rfloor$ for $i = 1, \ldots, k$. We also write $M$ for $\sum_{i=1}^k
m_i$.

The elements of the finite set $\{0, 1\}^M$ will be represented as
$\tau= (\tau_1, \ldots, \tau_k)$ where $\tau_i = (\tau_{i1}, \ldots,
\tau_{im_i}) \in\{0, 1\}^{m_i}$. For each $\tau\in\{0, 1\}^M$, we
specify $\theta^{\tau} \in\M$ in the following way. For
$s_{i-1} + 1 \leq u \leq s_i$, the quantity $\theta^{\tau}_u$ is defined
as
\begin{eqnarray*}
&&\theta_{0,i} + \frac{\beta_n}{m_i} \sum_{v=1}^{m_i}
(v - \tau_{iv}) I \bigl\{(v-1)l + 1 \leq u - s_{i-1} \leq vl
\bigr\} \\
&&\qquad{}+ \beta_n I \{s_{i-1} + m_i l + 1 \leq u
\leq s_i \}.
\end{eqnarray*}
Because $\theta$ is constant on the set $\{u \dvtx s_{i-1}
+ 1 \leq u \leq s_i\}$ where it takes the value $\theta_{0, i}$, it follows
that $\ell_{\infty}(\theta^{\tau},\theta) \leq\beta_n$. This implies
that $\theta^{\tau} \in\nbr(\theta)$ for every $\tau$ as
$\beta_n^2 \leq R(n;\theta)$.

Also, because of the assumption $\min_{2 \leq i \leq k} (\theta_{0,i} -
\theta_{0,i-1}) \geq\beta_n$, it is evident that each $\theta^{\tau
}$ is
nondecreasing. We will apply Assouad's lemma to $\theta^{\tau}, \tau
\in\{0, 1\}^M$. For $\tau, \tau' \in\{0, 1\}^M$, we have
%
\begin{equation}
\label{oro} \ell^2\bigl(\theta^{\tau},
\theta^{\tau'}\bigr) = \frac{1}{n} \sum_{i=1}^k
\sum_{v =
1}^{m_i} \frac{l\beta_n^2}{m_i^2} I \bigl
\{\tau_{iv} \neq \tau'_{iv} \bigr\} =
\frac{l\beta_n^2}{n} \sum_{i=1}^k
\frac{\ham(\tau_i, \tau'_i)}{m_i^2}.
\end{equation}
Because
\[
m_i \leq\frac{n_i}{l} \leq\frac{c_2n}{kl} \qquad\qt{for each $1
\leq i \leq k$},
\]
we have
%
\begin{equation}
\label{a1s} \ell^2\bigl(\theta^{\tau},
\theta^{\tau'}\bigr) \geq\frac{k^2 l^3 \beta
_n^2}{c_2^2 n^3} \sum
_{i=1}^k \ham\bigl(\tau_i,
\tau'_i\bigr) = \frac{k^2 l^3 \beta_n^2}{c_2^2
n^3} \ham\bigl(\tau,
\tau'\bigr).
\end{equation}
Also, from \eqref{oro}, we get
%
\begin{equation}
\label{bas} \ell^2\bigl(\theta^\tau,
\theta^{\tau'}\bigr) \leq\frac{l \beta_n^2}{n
(\min_{1 \leq i
\leq k} m_i^2)} \qquad\qt{when $\ham\bigl(\tau,
\tau'\bigr) = 1$}.
\end{equation}
The quantity $\min_i m_i^2$ can be easily bounded from below by noting
that $n_i/l < m_i + 1 \leq2 m_i$ and that $n_i \geq c_1 n/k$. This
gives
%
\begin{equation}
\label{ntr} \min_{1 \leq i \leq k} m_i \geq
\frac{c_1 n }{2 k l}.
\end{equation}
Combining the above inequality with \eqref{bas}, we deduce
\[
\ell^2\bigl(\theta^{\tau}, \theta^{\tau'}\bigr) \leq
\frac{4 k^2 l^3 \beta
_n^2}{c_1^2
n^3} \qquad\qt{whenever $\ham\bigl(\tau, \tau'\bigr) = 1$}.
\]
This and Pinsker's inequality give
%
\begin{eqnarray}
\label{a2s} \|\P_{\theta^{\tau}} - \P_{\theta^{\tau'}}\|_{\mathrm{TV}}^2
&\leq&\frac
{1}{2} D (\P_{{\theta^\tau}}\|\P_{\theta^{\tau'}}) =
\frac{n}{4 \sigma^2} \ell^2\bigl(\theta^{\tau},
\theta^{\tau'}\bigr)
\nonumber
\\[-8pt]
\\[-8pt]
\nonumber
&\leq&\frac{k^2 l^3 \beta
_n^2}{c_1^2 n^2
\sigma^2}
\end{eqnarray}
whenever $\ham(\tau, \tau') = 1$.

Inequalities \eqref{a1s} and \eqref{a2s} in conjunction with
Assouad's lemma give
\[
\mix_n(\theta) \geq\frac{M k^2 l^3 \beta_n^2}{8 c_2^2 n^3} \biggl(1 - \frac{k \beta_n l^{3/2}}{c_1 n \sigma}
\biggr).
\]
Because of \eqref{ntr}, we get $M = \sum_i m_i \geq k \min_i m_i
\geq
c_1 n/(2l)$, and thus
%
\begin{equation}
\label{sog} \mix_n(\theta) \geq\frac{c_1 k^2 l^2 \beta_n^2}{16 c_2^2 n^2} \biggl(1 -
\frac{k \beta_n l^{3/2}}{c_1 n \sigma} \biggr).
\end{equation}
The value of the integer $l$ will now be specified. We take
%
\begin{equation}
\label{tar} l = \biggl( \frac{c_1 n \sigma}{ 2 k \beta_n} \biggr)^{2/3}.
\end{equation}
Because $\min_i n_i \geq c_1 n/k$, we can ensure that $1 \leq l \leq
\min_i n_i$ by requiring that
\[
1 \leq \biggl(\frac{c_1 n \sigma}{ 2 k \beta_n} \biggr)^{2/3} \leq \frac{c_1 n}{k}.
\]
This gives rise to two lower bounds for $n$ which are collected
in \eqref{tie.con}.

As a consequence of \eqref{tar}, we get that $l^{3/2} \leq c_1 n
\sigma/(2 k \beta_n)$, which ensures that the term inside the
parentheses on the right-hand side of \eqref{sog} is at least
$1/2$. This gives
%
\begin{equation}
\label{rar} \mix_n(\theta) \geq\frac{c_1 k^2 l^2 \beta_n^2}{32 c_2^2
n^2} \geq
\frac{c_1^{7/3}}{2^{19/3}
c_2^2} \frac{k \sigma^2}{n} \biggl(\log\frac{en}{k}
\biggr)^{1/3}.
\end{equation}
To complete the proof, we use Theorem~\ref{vt}. Specifically, the
second inequality in~\eqref{vt.eq} gives
\[
\frac{k \sigma^2}{n} \geq\frac{R(n;\theta)}{16} \biggl(\log \frac{en}{k}
\biggr)^{-1}.
\]
The proof is complete by combining the above inequality
with \eqref{rar}.
\end{pf}

\section{Risk bound under model misspecification}\label{Misspec}
We consider isotonic regression under the misspecified setting where
the true sequence is not necessarily nondecreasing. Specifically,
consider model \eqref{eq:Mdl2} where now the true sequence $\theta$
is not necessarily assumed to be in $\M$. We study the behavior of the
LSE $\hat{\theta} = \hat{\theta}(Y; \M)$. The goal of this section is
to prove an inequality analogous to \eqref{pyaa} for model
misspecification. It turns out here that the LSE is really estimating
the nondecreasing projection of $\theta$ on $\M$ defined as
$\tilde{\theta} \in\M$ that minimizes $\ell^2(t, \theta)$ over $t
\in\M$. From~\cite{RWD88}, Chapter~1, it follows that
%
\begin{equation}
\label{eq:Misspec} \tilde\theta_j = \min_{l \geq j} \max
_{k \leq j} \bar{\theta}_{k, l} \qquad\qt{for $1 \leq j \leq n$},
\end{equation}
where $\bar\theta_{k,l}$ is as defined in \eqref{eq:ShortHand}.

We define another measure of variation for $t \in\M$ with respect to
an interval partition $\pi= (n_1, \ldots, n_k)$:
\[
S_{\pi} (t) = \Biggl( \frac{1}{n}\sum
_{i =
1}^{k} \sum_{j=s_{i-1} + 1}^{s_i}
(t_{s_i} - t_j )^2 \Biggr)^{1/2},
\]
where $s_0 = 0$ and $ s_i = n_1 + \cdots + n_i$ for $1 \leq i \leq
k$. It is easy to check that $S_{\pi} (t) \leq V_{\pi} (t)$ for
every $t \in\M$. The following is the main result of this
section. The proofs of all the results in this section can be found in
the supplementary material \cite{ChatEtAlSupp14}.

\begin{theorem}\label{thm:Misrafi}
For every $\theta\in\R^n$, the LSE satisfies
%
\begin{equation}
\label{miss.eq} \E_{\theta} \ell^2(\tilde{\theta}, \hat{\theta})
\leq4 \inf_{\pi\in{\Pi}} \biggl(S^2_{\pi} (
\tilde{\theta}) + \frac{4 \sigma k(\pi)}{n} \log \frac{en}{k(\pi)} \biggr) \leq R(n;
\tilde{\theta}).\vadjust{\goodbreak}
\end{equation}
\end{theorem}
%
\begin{remark}
By Theorem~\ref{sabya2}, the quantity $R(n; \tilde{\theta})$ is
bounded from above by $(\sigma^2 V(\tilde{\theta})/n)^{2/3}$ up to a
logarithmic multiplicative factor in $n$. Therefore,
Theorem~\ref{thm:Misrafi} implies that the LSE $\hat{\theta}$
converges to the projection of $\theta$ onto the space of monotone
vectors at least the $n^{-2/3}$ rate, up to a logarithmic factor in
$n$. The convergence rate will be much faster if $k(\tilde{\theta})$
is small or if $\tilde{\theta}$ is well approximated by a monotone
vector $\alpha$ with small $k(\alpha)$.
\end{remark}

By taking $\pi$ in the infimum in the upper bound of \eqref{miss.eq}
to be the interval partition generated by $\tilde{\theta}$, we obtain
the following result which is the analogue of \eqref{amon} for model
misspecification.

\begin{corollary}
For every arbitrary sequence $\theta$ of length $n$ (not necessarily
nondecreasing),
\[
\E_{\theta} \ell^2(\tilde{\theta},\hat{\theta}) \leq
\frac{16
\sigma^2
k(\tilde{\theta})}{n} \log \frac{en}{k(\tilde{\theta})}.
\]
\end{corollary}

In the next pair of results, we prove two upper bounds on
$k(\tilde{\theta})$. The first result shows that $k(\tilde{\theta}) =
1$ (i.e., $\tilde{\theta}$ is constant) when $\theta$ is
nonincreasing, that is, $\theta_1 \ge\theta_2 \ge\cdots\ge
\theta_n$. This implies that the LSE converges to $\tilde{\theta}$ at
the rate $\sigma^2 \log(en)/n$ when $\theta$ is nonincreasing.

\begin{lemma}\label{appr}
$k(\tilde{\theta}) = 1$ if $\theta$ is nonincreasing.
\end{lemma}

To state our next result, let
%
\begin{equation}
\label{eq:block} b(t) := \sum_{i =1}^{n-1} I
\{t_i \neq t_{i+1} \} + 1 \qquad\qt{for $t \in\R^n$}.
\end{equation}
$b(t)$ can be interpreted as the number of constant blocks of
$t$. For example, when $n = 5$ and $t = (0,0,1,1,1,0)$,
$b(t) = 3$. Observe that $b(t) = k(t)$ for $t \in\M$.

\begin{lemma}\label{sabya3}
For any sequence $\theta\in\R^n$, we have $k(\tilde{\theta}) \leq
b(\theta)$.\vadjust{\goodbreak}
\end{lemma}

As a consequence of the above lemma, we obtain that for every $\theta
\in\R^n$, the quantity $\E_{\theta}
\ell^2(\hat{\theta}, \tilde{\theta})$ is bounded from above by
$(16b(\theta)\sigma^2/n) \log(en/b(\theta))$.

\section*{Acknowledgements} We thank the Associate Editor and the
anonymous referees
for constructive suggestions that significantly improved the paper.
\vspace*{6pt}

\begin{supplement}[id=suppA]
\stitle{Supplement to ``On risk bounds in isotonic and other shape
restricted regression problems''}
\slink[doi]{10.1214/15-AOS1324SUPP} 
\sdatatype{.pdf}
\sfilename{aos1324\_supp.pdf}
\sdescription{In the supplementary paper \cite{ChatEtAlSupp14} we provide
the proofs of Lemmas \ref{lili}, \ref{hol}, \ref{apps}, \ref{appr}
and \ref{sabya3} and
Theorems \ref{sabya2}, \ref{cst}, \ref{vt}, \ref{oval} and \ref
{thm:Misrafi}. We also state and prove Lemma~11.1,
which is used in the proof of Theorem~\ref{sabya2}, and Lemma~11.2, which is used
in the proof of Theorem~\ref{rafi}.}\vspace*{6pt}
\end{supplement}

%



\printaddresses
\end{document}